\documentclass{article}[11pt]
\usepackage{amsmath}
\usepackage{amsfonts}
\usepackage{amssymb}
\usepackage{color}

\newtheorem{theorem}{Theorem}
\newtheorem{lemma}[theorem]{Lemma}

\newtheorem{corollary}[theorem]{Corollary}

\newcommand{\Z}{{\mathbb Z}}

\newcommand{\R}{\mathbb R}

\def\IP{({\rm IP})_{N,\veb,\vel,\veu,f}}
\def\Graver{{\cal G}}

\def\Orthant_j{{\mathcal O}_{j}}

\usepackage{enumerate}
\def\ve#1{\mathchoice{\mbox{\boldmath$\displaystyle\bf#1$}}
{\mbox{\boldmath$\textstyle\bf#1$}}
{\mbox{\boldmath$\scriptstyle\bf#1$}}
{\mbox{\boldmath$\scriptscriptstyle\bf#1$}}}
\newcommand\vea{{\ve a}}
\newcommand\veb{{\ve b}}
\newcommand\vecc{{\ve c}}

\newcommand\veg{{\ve g}}

\newcommand\vel{{\ve l}}

\newcommand\veq{{\ve q}}

\newcommand\veu{{\ve u}}
\newcommand\vev{{\ve v}}
\newcommand\vew{{\ve w}}
\newcommand\vex{{\ve x}}
\newcommand\vey{{\ve y}}
\newcommand\vez{{\ve z}}

\newcommand{\boproof}{\textbf{Proof.} }
\newcommand{\eoproof}{\hspace*{\fill} $\square$ \vspace{5pt}}

\begin{document}
\setlength{\parindent}{0pt} \setlength{\parskip}{2ex plus 0.4ex
minus 0.4ex}

\title{A polynomial-time algorithm for optimizing over $N$-fold $4$-block decomposable integer programs}
\author{Raymond Hemmecke (Technische Universit\"at Munich, Germany) \and Matthias K\"oppe\thanks{Supported by NSF Grant}\ \ (University of California, Davis, USA)
  \and Robert Weismantel (Otto-von-Guericke-University, Magdeburg, Germany)}
\date{\today}

\maketitle

\begin{abstract}
In this paper we generalize $N$-fold integer programs and two-stage integer programs with $N$ scenarios to $N$-fold $4$-block decomposable integer programs. We show that for fixed blocks but variable $N$, these integer programs are polynomial-time solvable for any linear objective. Moreover, we present a polynomial-time computable optimality certificate for the case of fixed blocks, variable $N$ and any convex separable objective function. We conclude with two sample applications, stochastic integer programs with second-order dominance constraints and stochastic integer multi-commodity flows, which (for fixed blocks) can be solved in polynomial time in the number of scenarios and commodities and in the binary encoding length of the input data. In the proof of our main theorem we combine several non-trivial constructions from the theory of Graver bases. We are confident that our approach paves the way for further extensions.
\end{abstract}

\section{Introduction}

Let $A\in\Z^{d\times n}$ be a matrix. We associate with $A$ a finite set $\Graver(A)$ of vectors with remarkable properties. Consider the set $\ker(A)\cap\Z^n$. Then we put into $\Graver(A)$ all nonzero vectors $\vev\in\ker(A)\cap\Z^n$ that cannot be written as a sum $\vev=\vev'+\vev''$ of nonzero vectors $\vev',\vev''\in\ker(A)\cap\Z^n$ that lie in the same orthant (or equivalently, have the same sign pattern in $\{\geq\ve 0,\leq\ve 0\}^n$) as $\vev$. The set $\Graver(A)$ has been named the \emph{Graver basis} of $A$, since Graver \cite{Graver:75} introduced this set $\Graver(A)$ in $1975$ and showed that it constitutes an optimality certificate for a whole family of integer linear programs that share the same problem matrix, $A$. By this we mean, that $\Graver(A)$ provides an augmenting vector/step to any non-optimal feasible solution and hence allows the design of a simple augmentation algorithm to solve the integer linear program.

In the last $10$ years, a tremendous theoretical progress has been made in the theory of Graver bases. It has been shown that $\Graver(A)$ constitutes an optimality certificate for a much wider class of integer minimization problems, namely for those minimizing a concave or a separable convex objective function over $\{\vez:A\vez=\veb,\vel\leq\vez\leq\veu,\vez\in\Z^n\}$ \cite{DeLoera+Hemmecke+Onn+Rothblum+Weismantel,Murota+Saito+Weismantel,Onn+Rothblum}. Moreover, it has been shown that only polynomially many Graver basis augmentation steps are needed to find a feasible solution and to turn it into an optimal feasible solution \cite{Hemmecke:PSP,Hemmecke+Onn+Weismantel:oracle,Schulz+Weismantel}. Finally, based on the fundamental finiteness results for certain highly structured matrices $A$ ($N$-fold IPs and two- and multi-stage stochastic IPs) \cite{Aschenbrenner+Hemmecke,Hemmecke+Schultz,Hosten+Sullivant,Santos+Sturmfels}, it has been shown that concave and separable convex $N$-fold IPs and two- and multi-stage stochastic IPs can be solved in polynomial time \cite{DeLoera+Hemmecke+Onn+Weismantel,Hemmecke+Onn+Weismantel:oracle} for fixed blocks.

In this paper, we will combine the two cases of $N$-fold IPs and of two-stage stochastic IPs by considering problems with a problem matrix that is $N$-fold $4$-block decomposable as follows:
\[
\left(\begin{smallmatrix}C&D\\B&A\\\end{smallmatrix}\right)^{(N)}:=
\left(
\begin{array}{ccccc}
C & D & D & \cdots & D \\
B & A & 0 &   & 0 \\
B & 0 & A &   & 0 \\
\vdots &   &   & \ddots &   \\
B & 0 & 0 &   & A
\end{array}
\right)
\]
for some given $N\in\Z_+$ and $N$ copies of $A$. We call $\left(\begin{smallmatrix}C&D\\B&A\\\end{smallmatrix}\right)^{(N)}$ an $N$-fold $4$-block matrix. For $B=\ve 0$ and $C=\ve 0$ we recover the problem matrix of an $N$-fold IP and for $C=\ve 0$ and $D=\ve 0$ we recover the problem matrix of a two-stage stochastic IP.

Note that $N$-fold $4$-block decomposable matrices also arise in the context of combinatorial optimization \cite{Schrijver:86, Seymour}. More precisely, for totally unimodular matrices $C,A$ their $1$-sum is totally unimodular ($B=\ve 0, D=\ve 0$). Similarly, total unimodularity is preserved under the $2$-sum and $3$-sum composition.   Indeed, it can be verified that a repeated application of specialized $1$-sum, $2$-sum and $3$-sum compositions leads to a particular family of $N$-fold $4$-block decomposable matrices with structure regarding the matrices $B$ and $D$.

{\bf Example.} For matrices $C$ and $A$, column vector $\vea$ and row vector
$\veb^\intercal$ of appropriate dimensions, the $2$-sum of
$\left(\begin{smallmatrix}C&\vea\\\end{smallmatrix}\right)$ and
$\left(\begin{smallmatrix}\veb^\intercal\\A\\\end{smallmatrix}\right)$ gives
$\left(\begin{smallmatrix}C&\vea\veb^\intercal\\\ve
0&A\\\end{smallmatrix}\right)$. The $2$-sum of
$\left(\begin{smallmatrix}C&\vea\veb^\intercal & a\\\ve
0&A&0\\\end{smallmatrix}\right)$ and
$\left(\begin{smallmatrix}\veb^\intercal\\B\\\end{smallmatrix}\right)$ creates
the matrix
$\left(\begin{smallmatrix}C&\vea\veb^\intercal&\vea\veb^\intercal\\\ve 0&\ve
A&\ve 0\\\ve 0&\ve 0&A\\\end{smallmatrix}\right)$, which is the $2$-fold
$4$-block decomposable matrix
$\left(\begin{smallmatrix}C&\vea\veb^\intercal\\\ve
0&A\\\end{smallmatrix}\right)^{(2)}$. \eoproof

Our main result is the following.

\begin{theorem}\label{Theorem: ABCD}
Let $A\in\Z^{d_A\times n_A}$, $B\in\Z^{d_A\times n_B}$, $C\in\Z^{d_C\times n_B}$, $D\in\Z^{d_C\times n_A}$ be fixed matrices. For given $N\in\Z_+$ let $\vel,\veu\in\Z^{n_B+Nn_A}$, $\veb\in\Z^{d_C+Nd_A}$, and let $f:\R^{n_B+Nn_A}\rightarrow\R$ be a separable convex function and denote by $\hat f$ the maximum of $|f|$ over the feasible region of the convex integer minimization problem
\[
\IP:\qquad\min\left\{f(\vez):\left(\begin{smallmatrix}C&D\\B&A\\\end{smallmatrix}\right)^{(N)}\vez=\veb,\vel\leq\vez\leq\veu,\vez\in\Z^{n_B+Nn_A}\right\}.
\]
We assume that $f$ is given only by a comparison oracle that, when queried on $\vez$ and $\vez'$ decides whether $f(\vez)<f(\vez')$, $f(\vez)=f(\vez')$ or $f(\vez)>f(\vez')$. Then the following hold:
\begin{itemize}
  \item[(a)] There exists an algorithm that computes a feasible solution to $\IP$ or decides that no such solution exists and that runs in time polynomial in $N$, in the binary encoding lengths $\langle \vel,\veu,\veb\rangle$.
  \item[(b)] Given a feasible solution $\vez_0$ to $\IP$, there exists an algorithm that decides whether $\vez_0$ is optimal or finds a better feasible solution $\vez_1$ to $\IP$ with $f(\vez_1)<f(\vez_0)$ and that runs in time polynomial in $N$, in the binary encoding lengths $\langle \vel,\veu,\veb,\hat f\rangle$, and in the number of calls to the evaluation oracle for $f$.
  \item[(c)] If $f$ is linear, there exists an algorithm that finds an optimal solution to $\IP$ or decides that $\IP$ is infeasible or unbounded and that runs in time polynomial in $N$, in the binary encoding lengths $\langle \vel,\veu,\veb,\hat f\rangle$, and in the number of calls to the evaluation oracle for $f$.
\end{itemize}
\end{theorem}

This theorem generalizes a similar statement for $N$-fold integer programming and for two-stage stochastic integer programming. In these two special cases, one can even prove claim (c) of Theorem \ref{Theorem: ABCD} for all separable convex functions and for a certain class of separable convex functions, respectively. It is a fundamental open question, whether one can construct not only some augmenting vector for a given separable convex objective function $f$ in polynomially many steps but a best-improvement (or greedy) augmentation step $\alpha\vev$ with $\alpha\in\Z_+$ and $\vev\in\Graver\left(\left(\begin{smallmatrix}C&D\\B&A\\\end{smallmatrix}\right)^{(N)}\right)$. If this can be done, part (c) of Theorem \ref{Theorem: ABCD} can be extended from linear $f$ to a class of separable convex functions $f$ by applying the main result from \cite{Hemmecke+Onn+Weismantel:oracle}.

In fact, Theorem \ref{Theorem: ABCD} will be a consequence of the following structural result about $\Graver\left(\left(\begin{smallmatrix}C&D\\B&A\\\end{smallmatrix}\right)^{(N)}\right)$.

\begin{theorem}\label{Theorem: Polynomial bound for fixed ABCD on increase of Graver degree}
If $A\in\Z^{d_A\times n_A}$, $B\in\Z^{d_A\times n_B}$, $C\in\Z^{d_C\times n_B}$, $D\in\Z^{d_C\times n_A}$ are fixed matrices, then $\max\left\{\|\vev\|_1:\vev\in\Graver\left(\left(\begin{smallmatrix}C&D\\B&A\\\end{smallmatrix}\right)^{(N)}\right)\right\}$ is bounded by a polynomial in $N$.
\end{theorem}

In the next section, we present two applications of Theorem \ref{Theorem: ABCD}: stochastic integer programming with second-order dominance constraints and stochastic integer multi-commodity flows. While the first application has an $N$-fold $4$-block matrix as problem matrix, the second application can be modeled as an $N$-fold $4$-block IP after a suitable transformation. To state the result, we introduce the following type of matrices. For given $N\in\Z_+$ let
\[
\left[\begin{smallmatrix}A&B\\D&C\\\end{smallmatrix}\right]^{(N)}:=
\left(
\begin{array}{cccccc}
    A & & & B & \cdots & B \\
    & \ddots & & \vdots & & \vdots \\
    & & A & B & \cdots & B \\
    D & \cdots & D & C \\
    \vdots & & \vdots & & \ddots \\
    D & \cdots & D & & & C
\end{array}
\right),
\]
where we have $N$ copies of $A$ and of $C$. Then the following holds.

\begin{corollary}\label{Corollary: ABCD}
Let $A\in\Z^{d_A\times n_A}$, $B\in\Z^{d_A\times n_B}$, $C\in\Z^{d_C\times n_B}$, $D\in\Z^{d_C\times n_A}$ be fixed matrices. For given $N\in\Z_+$ let $\vel,\veu\in\Z^{N(n_A+n_B)}$, $\veb\in\Z^{N(d_A+d_C)}$, and let $f:\R^{N(n_A+n_B)}\rightarrow\R$ be a separable convex function and denote by $\hat f$ the maximum of $|f|$ over the feasible region of the convex integer minimization problem
\[
\IP':\qquad\min\left\{f(\vez):\left[\begin{smallmatrix}A&B\\D&C\\\end{smallmatrix}\right]^{(N)}\vez=\veb,\vel\leq\vez\leq\veu,\vez\in\Z^{N(n_A+n_B)}\right\}.
\]
We assume that $f$ is given only by a comparison oracle that, when queried on $\vez$ and $\vez'$ decides whether $f(\vez)<f(\vez')$, $f(\vez)=f(\vez')$ or $f(\vez)>f(\vez')$. Then the following hold:
\begin{itemize}
  \item[(a)] There exists an algorithm that computes a feasible solution to $\IP'$ or decides that no such solution exists and that runs in time polynomial in $N$, in the binary encoding lengths $\langle \vel,\veu,\veb\rangle$.
  \item[(b)] Given a feasible solution $\vez_0$ to $\IP'$, there exists an algorithm that decides whether $\vez_0$ is optimal or finds a better feasible solution $\vez_1$ to $\IP'$ with $f(\vez_1)<f(\vez_0)$ and that runs in time polynomial in $N$, in the binary encoding lengths $\langle \vel,\veu,\veb,\hat f\rangle$, and in the number of calls to the evaluation oracle for $f$.
  \item[(c)] If $f$ is linear, there exists an algorithm that finds an optimal solution to $\IP'$ or decides that $\IP'$ is infeasible or unbounded and that runs in time polynomial in $N$, in the binary encoding lengths $\langle \vel,\veu,\veb,\hat f\rangle$, and in the number of calls to the evaluation oracle for $f$.
\end{itemize}
\end{corollary}

We do now present problems to which Theorem \ref{Theorem: ABCD} and its Corollary \ref{Corollary: ABCD} apply. Thereafter, we prove our claims. Our proof of Theorem \ref{Theorem: ABCD} combines several non-trivial constructions from the theory of Graver bases. We are confident that our approach paves the way for further extensions.

\section{Applications}\label{Section: Applications}

In this section we present two $N$-fold $4$-block decomposable integer programming problems that are polynomial-time solvable for given fixed blocks and variable $N$ by Theorem \ref{Theorem: ABCD} and its Corollary \ref{Corollary: ABCD}.

\subsection{Stochastic integer multi-commodity flow}

Let there be $M$ integer (in contrast to continuous) commodities to be transported over a given network. While we assume that supply and demands are deterministic, we assume that the upper bounds for the capacities per edge are uncertain and given initially only via some probability distribution. The problem setup is as follows: first, we have to decide how to transport the $M$ commodities over the given network without knowing the true capacities per edge. Then, after observing the true capacities per edge, penalties have to be paid if the capacity is exceeded. Assuming that we have knowledge about the probability distributions of the uncertain upper bounds, we wish to minimize the costs for the integer multi-commodity flow plus the expected penalties to be paid for exceeding capacities. To solve this problem, we discretize as usual the probability distribution for the uncertain upper bounds into $N$ scenarios. Doing so, we obtain a (typically large-scale) (two-stage stochastic) integer programming problem with problem matrix
\[
\left(
\begin{array}{cccccccc}
    A & & & \ve 0 & \ve 0 &\cdots & \ve 0 & \ve 0\\
    & \ddots & & \vdots &\vdots & & \vdots &\vdots\\
    & & A & \ve 0 & \ve 0 & \cdots & \ve 0 & \ve 0 \\
    I & \cdots & I & I & -I&&&  \\
    \vdots & & \vdots & & & \ddots & & \\
    I & \cdots & I &   &   &   & I & -I
\end{array}
\right).
\]
Herein, $A$ is the node-edge incidence matrix of the given network, $I$ is an identity matrix of appropriate size, and the columns containing $-I$ correspond to the penalty variables. If the network is kept fix, $A$, $I$, and $-I$ are fix, too. As the problem matrix is simply $\left[\begin{smallmatrix}A&\ve 0\\I & \left(\begin{smallmatrix}I &-I\\\end{smallmatrix}\right)\\\end{smallmatrix}\right]^{(N)}$, we can apply Corollary \ref{Corollary: ABCD} and obtain the following.

\begin{theorem}
For given fixed network the two-stage stochastic integer linear multi-commodity flow problem is solvable in polynomial time in the number $M$ of commodities, in the number $N$ of scenarios, and in the encoding lengths of the input data.
\end{theorem}

\boproof The only issue that prevents us to apply Corollary \ref{Corollary: ABCD} directly is the fact that $M$ and $N$ are different. But by introducing additional commodities or scenarios, we can easily obtain an equivalent (bigger) problem with $M=N$ for which we can apply Corollary \ref{Corollary: ABCD}. If $M<N$, we introduce additional commodities with zero flow and if $M>N$, we take one scenario, copy it additional $M-N$ times and choose for each of these $M-N+1$ identical scenarios $1/(M-N+1)$ times the original cost vector. So, in total, these $M-N+1$ scenarios are equivalent to the one we started from. \eoproof

It should be noted that we can extend the problem and still get the same polynomiality result. For example, we may assume that we are allowed to change the routing of the $M$ commodities in the second-stage decision. Penalties could be enforced for the amount of change of the first-stage decision or only for the amount of additional flow on edges compared to the first-stage decision. Writing down the constraints and introducing suitable additional variables with zero lower and upper bounds, one obtains again a problem matrix that allows the application of Corollary \ref{Corollary: ABCD}.

\subsection{Stochastic integer programs with second-order dominance constraints}

Stochastic integer programs with second-order dominance constraints were introduced in \cite{Gollmer+Gotzes+Schultz}. Therein, in Proposition 3.1, the following mixed-integer linear program was obtained as a deterministic equivalent to solve the stochastic problem at hand. We refer the reader to \cite{Gollmer+Gotzes+Schultz} for the details.
\[
{\rm (SIP)}: \min\left\{\veg^\intercal\vex:
\begin{array}{rcl}
\vecc^\intercal\vex+\veq^\intercal\vey_{lk}-a_k & \leq & \vev_{lk}\;\;\;\forall l \forall k\\
T\vex+W\vey_{lk} & = & \vez_l\;\;\;\forall l \forall k\\
\sum_{l=1}^L \pi_l\vev_{lk} & \leq & \bar{a}_k\;\;\;\forall k\\
&&\vex\in X, \vey_{lk}\in\Z^{\bar m}_+\times\R^{m'}_+,\vev_{lk}\geq\ve 0\;\;\;\forall l \forall k
\end{array}
\right\}
\]
We assume now that all variables are integral and, for simplicity of exposition, we assume that the inequalities of the polyhedron $X$ are incorporated into the constraints $T\vex+W\vey_{lk}=\vez_l$. Moreover, we assume that all scenarios have the same probability, that is, $\pi_l=1/L$, $l=1,\ldots,L$.

\begin{theorem}
For given fixed matrices $T$ and $W$ and for fixed number $K$, problem {\rm (SIP)} is solvable in polynomial time in the number $L$ of (data) scenarios, and in the encoding lengths of the input data.
\end{theorem}

\boproof We transform the problem in such a way that Theorem \ref{Theorem: ABCD} can be applied.
First, we include the constraints $\vecc^\intercal\vex+\veq^\intercal\vey_{lk}-a_k\leq \vev_{lk}$ into the constraint $T\vex+W\vey_{lk}=\vez_l$ (by adding slack variables to get an equation). Then, we set $\bar T=\left(\begin{smallmatrix}T\\\vdots\\T\end{smallmatrix}\right)$ and $\bar W=\left(\begin{smallmatrix}W & &\\&\ddots&\\&&W\end{smallmatrix}\right)$, in which we use $K$ copies of $T$ and $W$, respectively. As $T$, $W$, and $K$ are assumed to be fixed, so are $\bar T$ and $\bar W$. With this, the problem matrix now becomes
\[
\left(
\begin{array}{ccccc}
    \ve 0 & I & \cdots & I & I \\
    \bar T & \bar W & & & \ve 0\\
    \vdots & & \ddots & & \vdots\\
    \bar T & & & \bar W & \ve 0\\
\end{array}
\right).
\]
Introducing suitable additional variables with zero lower and upper bounds, we obtain a problem matrix of the form $\left(\begin{smallmatrix}C&D\\B&A\\\end{smallmatrix}\right)^{(l)}$ with $A=\left(\begin{smallmatrix}\bar W &\ve 0\\\end{smallmatrix}\right)$, $B=\bar T$, $C=\ve 0$, and $D=\left(\begin{smallmatrix}I & I\\\end{smallmatrix}\right)$. Thus, we can apply Theorem \ref{Theorem: ABCD} and the result follows. \eoproof

\section{Proof of main results}

For a concise introduction to Graver bases (and to the results on $N$-fold IPs), including short proofs to the main results, we refer the reader to the survey paper by Onn \cite{Onn:survey}. In this section, we state and prove results on Graver bases needed for the proof of our main theorem in the next section. Let us start by bounding the $1$-norm of Graver basis elements of matrices with only one row. This lemma is a straight-forward consequence of Theorem 2 in \cite{Diaconis+Graham+Sturmfels}.

\begin{lemma}\label{Lemma: PPI degree bound}
Let $A\in\Z^{1\times n}$ be a matrix consisting of only one row and let $M$ be an upper bound on the absolute values of the entries of $A$. Then $\max\{\|\vev\|_1:\vev\in\Graver(A)\}\leq 2M-1$.
\end{lemma}

Let us now prove some more general degree bounds on Graver bases that we will use in the proof of the main theorem below.

\begin{lemma}\label{Lemma: General bound on increase of Graver degree}
Let $A\in\Z^{d\times n}$ and let $B\in\Z^{m\times n}$. Moreover, put $C:=\left(\begin{smallmatrix}A\\B\\\end{smallmatrix}\right)$. Then we have
\[
\max\{\|\vev\|_1:\vev\in\Graver(C)\}\leq\max\{\|\ve\lambda\|_1:\ve\lambda\in \Graver(B\cdot\Graver(A))\}\cdot\max\{\|\vev\|_1:\vev\in\Graver(A)\}.
\]
\end{lemma}

\boproof
Let $\vev\in\Graver(C)$. Then $\vev\in\ker(A)$ implies that $\vev$ can be written as a nonnegative integer linear sign-compatible sum $\vev=\sum\lambda_i \veg_i$ using Graver basis vectors $\veg_i\in\Graver(A)$. Adding zero components if necessary, we can write $\vev=\Graver(A)\ve\lambda$. We now claim that $\vev\in\Graver(C)$ implies $\lambda\in\Graver(B\cdot\Graver(A))$ and the result follows.

First, observe that $\vev\in\ker(B)$ implies $B\vev=B\cdot(\Graver(A)\ve\lambda)=(B\cdot\Graver(A))\ve\lambda=\ve 0$ and thus, $\ve\lambda\in\ker(B\cdot\Graver(A))$. If $\ve\lambda\not\in\Graver(B\cdot\Graver(A))$, then it can be written as a sign-compatible sum $\ve\lambda=\ve\mu+\ve\nu$ with $\ve\mu,\ve\nu\in\ker(B\cdot\Graver(A))$. But then
\[
\vev=(\Graver(A)\ve\mu)+(\Graver(A)\ve\nu)
\]
gives a sign-compatible decomposition of $\vev$ into two vectors $\Graver(A)\ve\mu,\Graver(A)\ve\nu\in\ker(C)$, contradicting the minimality property of $\vev\in\Graver(C)$. Hence, $\lambda\in\Graver(B\cdot\Graver(A))$ and the result follows.
\eoproof

We will employ the following simple corollary.

\begin{corollary}\label{Corollary: Special bound on increase of Graver degree}
Let $A\in\Z^{d\times n}$ and let $a^\intercal\in\Z^n$ be a row vector. Moreover, put $C:=\left(\begin{smallmatrix}A\\a\\\end{smallmatrix}\right)$. Then we have
\[
\max\{\|\vev\|_1:\vev\in\Graver(C)\}\leq\left(2\cdot\max\left\{|a^\intercal\vev|:\vev\in\Graver(A)\right\}-1\right)\cdot\max\{\|\vev\|_1:\vev\in\Graver(A)\}.
\]
In particular, if $M:=\max\{|a^{(i)}|:i=1,\ldots,n\}$ then
\[
\max\{\|\vev\|_1:\vev\in\Graver(C)\}\leq 2nM\left(\max\{\|\vev\|_1:\vev\in\Graver(A)\}\right)^2.
\]
\end{corollary}

\boproof
By Lemma \ref{Lemma: General bound on increase of Graver degree}, we already get
\[
\max\{\|\vev\|_1:\vev\in\Graver(C)\}\leq\max\{\|\ve\lambda\|_1:\ve\lambda\in \Graver(a^\intercal\cdot\Graver(A))\}\cdot\max\{\|\vev\|_1:\vev\in\Graver(A)\}.
\]
Now, observe that $a^\intercal\cdot\Graver(A))$ is a $1\times|\Graver(A)|$-matrix. Thus, the degree bound of primitive partition identities, Lemma \ref{Lemma: PPI degree bound}, applies, which gives
\[
\max\{\|\ve\lambda\|_1:\ve\lambda\in \Graver(a\cdot\Graver(A))\}\leq 2\cdot\max\left\{|a^\intercal\vev\|:\vev\in\Graver(A)\right\}-1,
\]
and thus, the first claim is proved. The second claim is a trivial consequence of the first. \eoproof

Let us now extend this corollary to a form that we need to prove Theorem \ref{Theorem: ABCD}.

\begin{corollary}\label{Corollary: Recursive bound on increase of Graver degree}
Let $A\in\Z^{d\times n}$ and let $B\in\Z^{m\times n}$. Let the entries of $B$ be bounded by $M$ in absolute value. Moreover, put $C:=\left(\begin{smallmatrix}A\\B\\\end{smallmatrix}\right)$. Then we have
\[
  \max\{\|\vev\|_1:\vev\in\Graver(C)\}\leq (2nM)^{2^m-1}\left(\max\{\|\vev\|_1:\vev\in\Graver(A)\}\right)^{2^m}.
\]
\end{corollary}

\boproof This claim follows by simple induction, adding one row of $B$ at a time, and by using the second inequality of Corollary \ref{Corollary: Special bound on increase of Graver degree} to bound the sizes of the intermediate Graver bases in comparison to the Graver basis of the matrix with one row of $B$ less. \eoproof

We are now ready to prove Theorem \ref{Theorem: ABCD}. We start out by considering the submatrix $\left(\begin{smallmatrix}\ve 0& \ve 0\\B&A\\\end{smallmatrix}\right)^{(N)}$. A main result from \cite{Hemmecke+Schultz} is the following.

\begin{lemma}\label{Lemma: Finite bound for SIP Graver bases}
  Let $A\in\Z^{d_A\times n_A}$ and $B\in\Z^{d_A\times n_B}$. There exists a number $g\in\Z_+$ depending only on $A$ and $B$ but not on $N$ such that for every $N\in\Z_+$ and for every $\vev\in\Graver\left(\left(\begin{smallmatrix}\ve 0&\ve 0\\B&A\\\end{smallmatrix}\right)^{(N)}\right)$, the components of $\vev$ are bounded by $g$ in absolute value. In particular, $\|\vev\|_1\leq (n_B+Nn_A)g$ for all $\vev\in\Graver\left(\left(\begin{smallmatrix}\ve 0&\ve 0\\B&A\\\end{smallmatrix}\right)^{(N)}\right)$.
\end{lemma}

Combining this result with Corollary \ref{Corollary: Recursive bound on increase of Graver degree}, we get a bound for the $1$-norms of the Graver basis elements of $\left(\begin{smallmatrix}C&D\\B&A\\\end{smallmatrix}\right)^{(N)}$.

\begin{corollary}\label{Corollary: Polynomial bound for fixed ABCD on increase of Graver degree}
Let $A\in\Z^{d_A\times n_A}$, $B\in\Z^{d_A\times n_B}$, $C\in\Z^{d_C\times n_B}$, $D\in\Z^{d_C\times n_A}$ be given matrices. Moreover, let $M$ be a bound on the absolute values of the entries in $C$ and $D$, and let $g\in\Z_+$ be the number from Lemma \ref{Lemma: Finite bound for SIP Graver bases}. Then for any $N\in\Z_+$ we have
\begin{eqnarray*}
  \max\left\{\|\vev\|_1:\vev\in\Graver\left(\left(\begin{smallmatrix}\ve 0&\ve 0\\B&A\\\end{smallmatrix}\right)^{(N)}\right)\right\} & \leq & (2(n_B+Nn_A)M)^{2^{d_C}-1}\left(\max\left\{\|\vev\|_1:\vev\in\Graver\left(\left(\begin{smallmatrix}\ve 0&\ve 0\\B&A\\\end{smallmatrix}\right)^{(N)}\right)\right\}\right)^{2^{d_C}}\\
& \leq & (2(n_B+Nn_A)M)^{2^{d_C}-1}\left((n_B+Nn_A)g\right)^{2^{d_C}}.
\end{eqnarray*}

If $A$, $B$, $C$, $D$ are fixed matrices, then $\max\left\{\|\vev\|_1:\vev\in\Graver\left(\left(\begin{smallmatrix}C&D\\B&A\\\end{smallmatrix}\right)^{(N)}\right)\right\}$ is bounded by a polynomial in $N$.
\end{corollary}

\boproof
While the first claim is a direct consequence of Lemma \ref{Lemma: Finite bound for SIP Graver bases} and Corollary \ref{Corollary: Recursive bound on increase of Graver degree}, the polynomial bound for fixed matrices $A$, $B$, $C$, $D$ and varying $N$ follows immediately by observing that $n_A,n_B,d_C,M,g$ are constants as they depend only on the fixed matrices $A$, $B$, $C$, $D$.
\eoproof

Note that the second claim of Corollary \ref{Corollary: Polynomial bound for fixed ABCD on increase of Graver degree} is exactly Theorem \ref{Theorem: Polynomial bound for fixed ABCD on increase of Graver degree}. Now we are ready to prove our main theorem.

{\bf Proof of Theorem \ref{Theorem: ABCD}.}
Let $N\in\Z_+$, $\vel,\veu\in\Z^{n_B+Nn_A}$, $\veb\in\Z^{d_C+Nd_A}$, and a separable convex function $f:\R^{n_B+Nn_A}\rightarrow\R$ be given. To prove claim (a), observe that one can turn any integer solution to $\left(\begin{smallmatrix}C&D\\B&A\\\end{smallmatrix}\right)^{(N)}\vez=\veb$ (which can be found in polynomial time using for example the Hermite normal form of $\left(\begin{smallmatrix}C&D\\B&A\\\end{smallmatrix}\right)^{(N)}$) into a feasible solution (that in addition fulfills $\vel\leq\vez\leq\veu$) by a sequence of linear integer programs (with the same problem matrix $\left(\begin{smallmatrix}C&D\\B&A\\\end{smallmatrix}\right)^{(N)}$) that ``move'' the components of $\vez$ into the direction of the given bounds, see \cite{Hemmecke:PSP}. This step is similar to phase I of the Simplex Method in linear programming. In order to solve these linear integer programs, it suffices (by the result of \cite{Schulz+Weismantel}) to find Graver basis augmentation vectors from $\Graver\left(\left(\begin{smallmatrix}C&D\\B&A\\\end{smallmatrix}\right)^{(N)}\right)$ for a directed augmentation oracle. So, claim (b) will imply both claim (a) and claim (c).

Let us now assume that we are given a feasible solution $\vez_0=(\vex,\vey_1,\ldots,\vey_N)$ and that we wish to decide whether there exists another feasible solution $\vez_1$ with $f(\vez_1)<f(\vez_0)$. By the main result in \cite{Murota+Saito+Weismantel}, it suffices to decide whether there exists some vector $\vev=(\hat\vex,\hat\vey_1,\ldots,\hat\vey_N)$ in the Graver basis of $\left(\begin{smallmatrix}C&D\\B&A\\\end{smallmatrix}\right)^{(N)}$ such that $\vez_0+\vev$ is feasible and $f(\vez_0+\vev)<f(\vez_0)$. By Corollary \ref{Corollary: Polynomial bound for fixed ABCD on increase of Graver degree} and by the fact that $n_B$ is constant, there is only a polynomial number of candidates for the $\hat\vex$-part of $\vev$. For each such candidate $\hat\vex$, we can find a best possible choice for $\hat\vey_1,\ldots,\hat\vey_N$ by solving the following separable convex $N$-fold IP:
\[
\min\left\{f\left(\left(\begin{smallmatrix}\vex+\hat\vex\\\vey_1+\hat\vey_1\\\vdots\\\vey_N+\hat\vey_N\end{smallmatrix}\right)\right):
\left(\begin{smallmatrix}C&D\\B&A\\\end{smallmatrix}\right)^{(N)}
\left(\begin{smallmatrix}\vex+\hat\vex\\\vey_1+\hat\vey_1\\\vdots\\\vey_N+\hat\vey_N\end{smallmatrix}\right)=
\veb,\vel\leq\left(\begin{smallmatrix}\vex+\hat\vex\\\vey_1+\hat\vey_1\\\vdots\\\vey_N+\hat\vey_N\end{smallmatrix}\right)\leq\veu, \vey_1,\ldots,\vey_N\in\Z^{n_A}\right\},
\]
for given $\vez_0=(\vex,\vey_1,\ldots,\vey_N)$ and $\hat\vex$. Observe that the problem $\IP$ does indeed simplify to a separable convex $N$-fold IP with problem matrix $\left(\begin{smallmatrix}0&D\\0&A\\\end{smallmatrix}\right)^{(N)}$ because $\vez_0=(\vex,\vey_1,\ldots,\vey_N)$ and $\hat\vex$ are fixed. For fixed matrices $A$ and $D$, however, each such $N$-fold IP is solvable in polynomial time \cite{Hemmecke+Onn+Weismantel:oracle}. If the $N$-fold IP is feasible and if for the resulting optimal vector $\vev:=(\hat\vex,\hat\vey_1,\ldots,\hat\vey_N)$ we have $f(\vez_0+\vev)\geq f(\vez_0)$, then no augmenting vector can be constructed using this particular choice of $\hat\vex$. If on the other hand we have $f(\vez_0+\vev)< f(\vez_0)$, then $\vev$ is a desired augmenting vector for $\vez_0$ and we can stop. As we solve polynomially many polynomially solvable $N$-fold IPs, claim (b) and thus also claims (a) and (c) follow. \eoproof

{\bf Proof of Corollary \ref{Corollary: ABCD}.} To prove Corollary \ref{Corollary: ABCD}, observe that after introducing additional variables, problem $\IP'$ can be modeled as an $N$-fold $4$-block IP and is thus polynomial-time solvable by Theorem \ref{Theorem: ABCD}. First, write the constraint $\left[\begin{smallmatrix}A&B\\D&C\\\end{smallmatrix}\right]^{(N)}\vez=\veb$ in $\IP'$ as follows:
\[
\left(
\begin{array}{cccccc}
    A & & & B & \cdots & B \\
    & \ddots & & \vdots & & \vdots \\
    & & A & B & \cdots & B \\
    D & \cdots & D & C \\
    \vdots & & \vdots & & \ddots \\
    D & \cdots & D & & & C
\end{array}
\right)\left(\begin{smallmatrix}\vex_1\\\vdots\\\vex_N\\\vey_1\\\vdots\\\vey_N\\
\end{smallmatrix}\right)=\left(\begin{smallmatrix}\veb_1\\\vdots\\\veb_N\\\veb_{N+1}\\\vdots\\\veb_{2N}\\
\end{smallmatrix}\right).
\]
Now introduce variables $\vew_x=\sum_{i=1}^N \vex_i$ and $\vew_y=\sum_{i=1}^N \vey_i$. Then we get the new constraints
\[
\left(
\begin{array}{cc|ccccccc}
    -I & & I & & I & & \cdots & I \\
    & -I & & I & & I & \cdots & & I \\
\hline
    D & &   & C \\
    & B & A &  \\
    D & & & & & C \\
    & B & & & A \\
    \vdots & & & & & & \ddots \\
    D & & & & & & & & C \\
    & B & & & & & & A
\end{array}
\right)\left(\begin{smallmatrix}\vew_x\\\vew_y\\\vex_1\\\vey_1\\\vdots\\\vex_N\\\vey_N\\\end{smallmatrix}\right)=
\left(\begin{smallmatrix}\ve 0\\\ve 0\\\veb_1\\\vdots\\\veb_N\\\veb_{N+1}\\\vdots\\\veb_{2N}\\
\end{smallmatrix}\right).
\]
Hence, $\IP'$ can be modeled as an $N$-fold $4$-block decomposable IP and thus, Corollary \ref{Corollary: ABCD} follows by applying Theorem \ref{Theorem: ABCD} to this transformed integer program. \eoproof

{\bf Acknowledgments.} We wish to thank R\"udiger Schultz for valuable comments on Section \ref{Section: Applications} and for pointing us to \cite{Gollmer+Gotzes+Schultz}.


\end{document}